\documentclass[12pt]{article}
\usepackage{amsmath, amsthm}
\usepackage{amsfonts,amssymb}
\usepackage{epsfig}
\usepackage{graphicx}
\usepackage{subfigure}
\newcommand{\R}{{\mathbb R}}
\newcommand{\bb}{\mathbf}

\newcommand{\m}{\medskip}
\newtheorem{theorem}{Theorem}

\newtheorem{lemma}{Lemma}

\theoremstyle{definition}

\newtheorem{example}{Example}

\def\sqr#1#2{{\vcenter{\vbox{\hrule height.#2pt
 \hbox{\vrule width.#2pt height#1pt \kern#1pt
 \vrule width.#2pt}
 \hrule height.#2pt}}}}
\def\square{\mathchoice\sqr68\sqr68\sqr{2.1}3\sqr{1.5}3}
\date{ }

\title{Fractal Homeomorphism for Bi-affine Iterated Function Systems} 
\author{
Michael Barnsley \\  The
Australian National University \\Canberra, Australia \\ \tt michael.barnsley@anu.edu.au
\\ \\
Andrew Vince \\ Department of Mathematics, University
of Florida \\ Gainesville, FL, USA \\ \tt avince@ufl.edu }

\begin{document}
\maketitle

\begin{abstract}  The paper concerns fractal homeomorphism between the attractors of two bi-affine iterated function systems.  After a general discussion of bi-affine functions, conditions are provided under which a bi-affine iterated function system is contractive, thus guaranteeing an attractor.  After a general discussion of fractal homeomorphism, fractal homeomorphisms are constructed for a specific type of bi-affine iterated function system.   
\end{abstract}

\section{Introduction} The purpose of this paper is to investigate bi-affine iterated function systems and  fractal homeomorphisms between their attractors.   The class of  bi-affine functions from $\R^2$ to $\R^2$ is more general than affine transformations but less general than quadratic transformations. These are the functions $f \, : \, \R^2 \rightarrow \R^2$  that are, for a fixed $x$ or a fixed $y$, affine in the other variable:
\begin{equation} \label{proportion}  \begin{aligned}  f((1-\alpha)x_1+ \alpha x_2,y) &= (1-\alpha)\, f(x_1,y) + \alpha \, f(x_2,y) \quad \text{  and   } \\ f(x, (1-\alpha)y_1+ \alpha y_2,) &= (1-\alpha)\, f(x,y_1) + \alpha \, f(x,y_2)  \end{aligned} 
\end{equation}
for all $x_1,x_2,y_1,y_2, \alpha \in \R$.  Interpreted geometrically, these equations mean that 
\begin{enumerate}
\item Horizontal and vertical lines are taken to lines, and
\item proportions along horizontal and vertical lines are preserved.
\end{enumerate}
This elementary class of functions, with connections to classic geometric results of Brianchon and Lampert dating back to the 18th century, proves extremely versatile for the applications described in this paper.  

\begin{figure}[htb] \label{lena}
\begin{center}
\includegraphics[width=4.5cm, keepaspectratio]{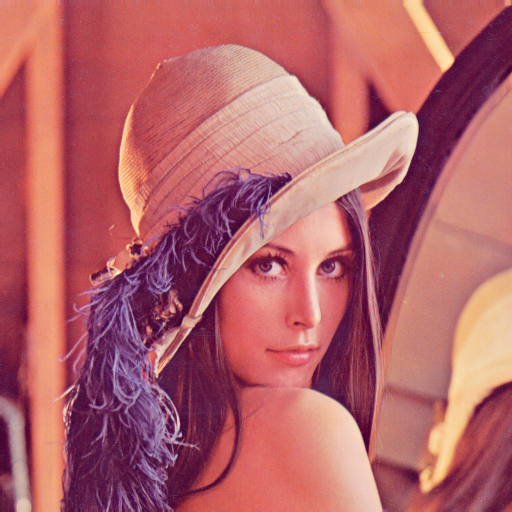} \hskip 8mm \includegraphics[width=4.5cm, keepaspectratio]{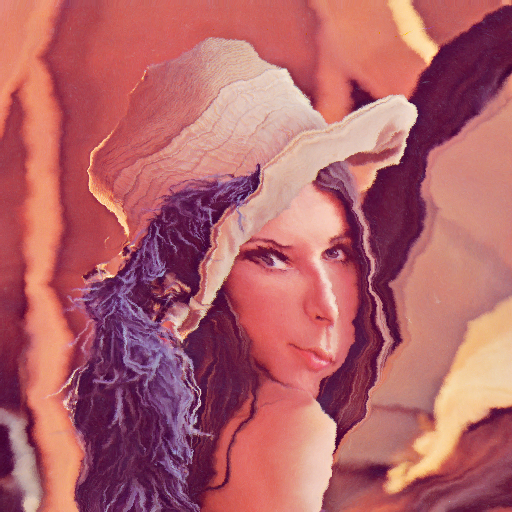}  \hskip 8mm
\includegraphics[width=4.5cm, keepaspectratio]{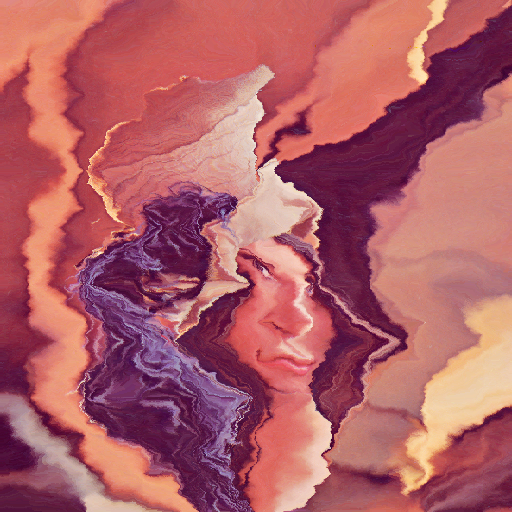} 
\caption{Two fractal homeomorphisms applied to the original picture.}
\end{center} 
\end{figure}

Our main motivation for investigating bi-affine functions comes from the representation and transformation of certain fractal images.  A standard method for constructing a deterministic self-referential fractal is by an iterated function system (IFS).  The attractor of the IFS is usually a fractal.  Barnsley \cite{ B1, B2} has introduced a method for transforming the attractor of one IFS to the attractor of another IFS, a method that has applications to digital imaging such as image encryption, filtering, compression, watermarking, and various special effects.  Figure~\ref{lena}, explained in more detail in Section~\ref{frSection}, is obtained by applying such transformations - called fractal homeomorphisms. In constructing
 a fractal homeomorphism it is convenient to use an IFS whose maps have nice geometric properties but that are not too complicated. Linear transformations are easy to work with and have the property that lines are taken to lines.  Affine transformations are not much more complicated than linear transformations and do not have the restriction that the origin be taken to the origin. There is a tradeoff; the more properties required, the more complicated the mapping. In the fractal geometry literature, in particular for fractals constructed using an iterated function system, affine transformations are frequently used.  For the applications described in Section~\ref{frSection}, requiring slightly more general functions,  bi-affine functions are ideally suited. \m

The paper is organized as follows.  The geometry of bi-affine functions is the subject of Section~\ref{BF}. 
Theorem~\ref{properties} gives basic properties of a bi-affine function, in particular properties of the folding line and
folding parabola.  Theorem~\ref{image} gives a geometric construction for finding  the image of a given point under a bi-affine 
function and  precisely describes the 2-to-1 nature of a bi-affine function.  
 Section~\ref{IFS}  provides background on iterated function systems and their attractors.  It is a classic result  that,
if an IFS  is contractive, then it has an attractor.   Theorem~\ref{contraction} gives fairly general conditions under which a bi-affine IFS is contractive.   Fractal homeomorphism between the attractors to two IFSs is the subject of Section~\ref{FT}.   The construction of a fractal homeomorphism depends on finding shift invariant sections of the coding maps of the two IFSs.  Theorem~\ref{mask} states that any such shift invariant section comes from a mask.  (The terms coding map, section, shift invariant, and mask are defined in Section~\ref{FT}.)   Theorem~\ref{pi-tau} concerns how to obtain a fractal homeomorphism between two attractors from the respective sections.  The visual representation of a fractal homeomorphism for bi-affine IFSs is the subject of  Section~\ref{frSection}.  Theorem~\ref{picture} states that a  particular type of bi-affine IFS can be used to construct fractal homeomorphisms.  The main issue is proving the continuity of the map and its inverse. The pictures in Figure~\ref{lena} are obtained by this method.  

\section{Geometry of Bi-affine Functions} \label{BF}

Boldface letters represent vectors in $\R^2$.  A function $f \, : \, \R^2 \rightarrow \R^2$ is called {\it bi-affine} if it has the form 
\begin{equation} \label{bi-affine} f(x,y) = {\bb a} + {\bb b}x + {\bb c} y + {\bb d} xy. \end{equation}
It is easy to verify that the class of bi-affine functions is exactly the class characterized by the two properties listed in the introduction.  The class of bi-affine functions is not closed under composition, but the composition of a bi-affine and an affine function is bi-affine. 

 Call a bi-affine function {\it non-degenerate} if $\bb d \neq \bb 0$ and neither $\bb b$ nor $\bb c$ is a scalar multiple of ${\bb d}$. In particular, neither $\bb b$ nor $\bb c$ is the zero vector. If $\bb d=0$, then $f$ is ``degenerate" in the sense that it is affine and well understood.  If $\bb d \neq \bb 0$ and both $\bb b$ and $\bb c$ are scalar multiples of ${\bb d}$, then $f$ is ``degenerate" in the sense that the range of $f$ degenerates to a line.  If $\bb d \neq \bb 0$ and just one of $\bb b$ and $\bb c$ is a scalar multiple of ${\bb d}$, then $f$ is ``degenerate" in the sense that the image of the folding line, as defined in the next section, is just a point, a fact that can be verified by equation~\ref{par} in the proof of Theorem~\ref{properties}.

 Basic properties of bi-affine functions are described in this section.  According to statement 2 of Theorem~\ref{properties} below, the image of a line $L$ under a bi-affine function  is a parabola. Such a parabola can be {\it degenerate} in the sense that it is either a line (focal distance $0$) or a line that doubles back on itself (focal distance $\infty$).  The first case occurs if and only if $L$ is parallel to either the $x$ or $y$-axis.  An example of the second case is
$(X,Y) = f(x,y) = (- x + xy, -1-y+xy)$, in which case the image of the line $y=x+1$ is given by the parametric equation $X = t^2, \, Y = t^2-2$, a degenerate parabola with vertex at $(0,-2)$ that doubles back along the line $y=x-2$.  

For vectors $\bb a$ and $\bb b$, let $|\bb a \; \bb b|$ denote the determinant
of the matrix whose columns are $\bb a$ and $\bb b$.  For a non-degenerate bi-affine functon $f$, call the line $L_f$ with equation
\begin{equation} \label{folding} |\bb b \; \bb d| \,x+ |\bb d \; \bb c| \, y = |\bb c \; \bb b| \end{equation} the {\it folding line} of the function $f$.  Note that, by the non-degeneracy of $f$, neither $|\bb b \; \bb d|$ nor $|\bb d \; \bb c|$ is zero, and hence the folding line is not parallel to either the $x$ or $y$-axis. The terminology ``folding line" is justified by the next theorem, in particular statement 4.   The image $P_f = f(L_f)$ of the folding line is called the {\it folding parabola}, which, according to statement 3 is non-degenerate. Let $L_{+}$ and $L_{-}$ denote the closed half spaces above and below $L_f$, respectively.  

\begin{theorem} \label{properties} If $f$ is a non-degenerate bi-affine map, then 
\begin{enumerate}
\item  any line parallel to either coordinate axis is mapped to a line;   
\item  any line is mapped to a, possibly degenerate, parabola;  
\item  the folding parabola $P_f := f(L_f)$ is non-degenerate;
\item  the map $f$ is injective when restricted to either $L_{+}$ or $L_{-}$. 
\end{enumerate}
\end{theorem}

\proof  Statement (1) follows from the fact that, for a fixed $x$ or a fixed $y$, a bi-affine function is an affine function in the other variable.  

The image of a line under a bi-affine transformation has a parametric equation of the form $f(t) = 
\bb u + \bb v \,t + \bb w \,t^2$, which is a, possibly degenerate, parabola.

To show that $P_f$ is non-degenerate, first note that the image of the line $L_f$ under $f$ is given by the parametric equation $(X(x), Y(x))$ with parameter $x$ by 
\begin{equation} \label{par} \begin{aligned} (X,Y) = &  f(x,x) = \frac{1}{ |\bb d \, \bb c|} \left [(|\bb d \, \bb c|\, \bb a + |\bb c \, \bb b| \, \bb c)+
(|\bb d \, \bb c|\, \bb b - |\bb b \, \bb d|\, \bb c +|\bb c \, \bb b|\, \bb d)\, x + |\bb b \, \bb d|\,\bb d \, x^2 \right ]\\
= &  \frac{1}{ |\bb d \, \bb c|} \left [(|\bb d \, \bb c|\, \bb a + |\bb c \, \bb b| \, \bb c)+
2 |\bb d \, \bb b| \,\bb c \, x+ |\bb b \, \bb d|\,\bb d \, x^2 \right ],\end{aligned} \end{equation}
the first equality simply by substituting from equation~\ref{folding} into equation~\ref{bi-affine} and the second equality by a direct calculation. The tangent vector to this parabola is 
$$T(x) = \frac{2\,|\bb d \, \bb b|}{ |\bb d \, \bb c|} ( \bb c -  \bb d \, x ).$$
By the non-degeneracy of $f$, the vectors $\bb c$ and $\bb d$ are linearly independent, implying that
the the direction of the tangent vector is not constant. Hence the parabola is not degenerate.

The Jacobian determinant of a bi-affine function is $|\bb b \; \bb d| \,x+ |\bb d \; \bb c| \, y - |\bb c \; \bb b|$, which is nonzero except on the folding line.  It follows from the inverse function theorem that $f$ is injective when restricted to either $L_{+}$ or $L_{-}$.
\qed \m

A parabola $P$ divides the plane into two regions; let $\widehat P$ denote the closed region ``outside" P (and including $P$). The set $\widehat{P_f}$ will be called simply the {\it parabolic region of} $f$. For a point $\bb p = (x,y) \in \R^2$, let 
\begin{equation} {\bb p} ^* = \left ( \frac{|\bb c \; \bb b| - |\bb d \; \bb c| \,y}{|\bb b \; \bb d|} \, , \, \frac{|\bb c \; \bb b|-|\bb b \; \bb d| \,x}{|\bb d \; \bb c|} \right ).\end{equation}  

\begin{figure}[htb] 
\vskip -2mm
\begin{center}
\includegraphics[width=1.5in, keepaspectratio] {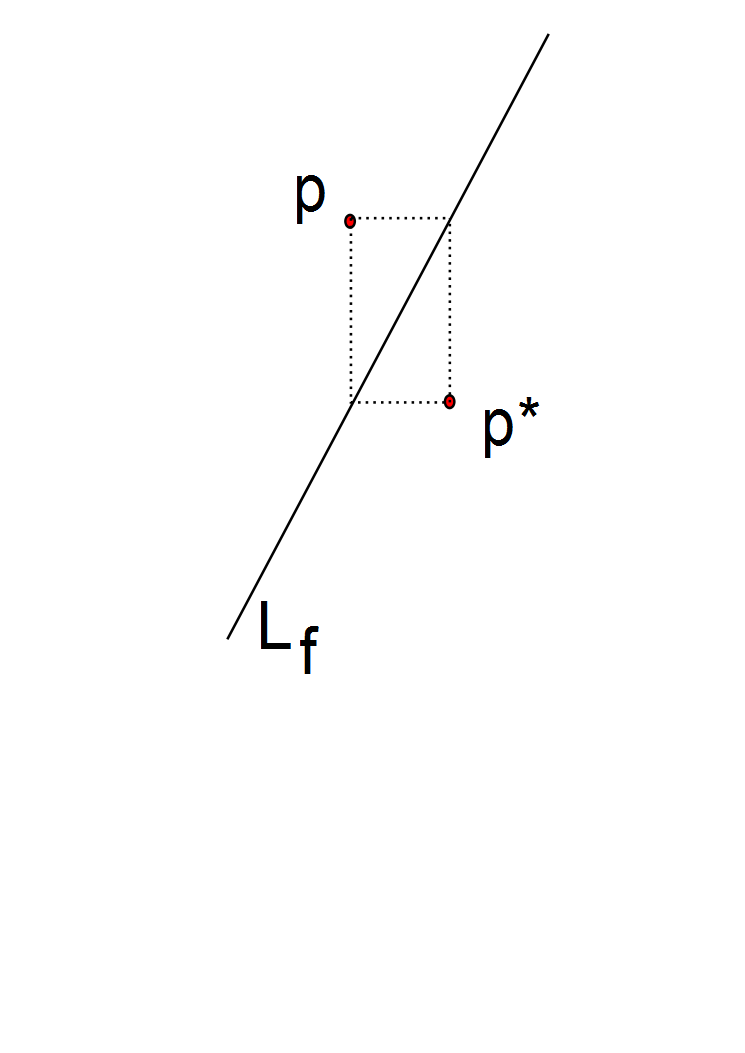}
\vskip -2cm
\caption{Folding Line $L_f:  f({\bb p}) ={\bb p^*}$ .}
\label{2-1}
\end{center}
\end{figure}

\begin{figure}[htb]  
\vskip -4mm
\begin{center}
\includegraphics[width=3in, keepaspectratio] {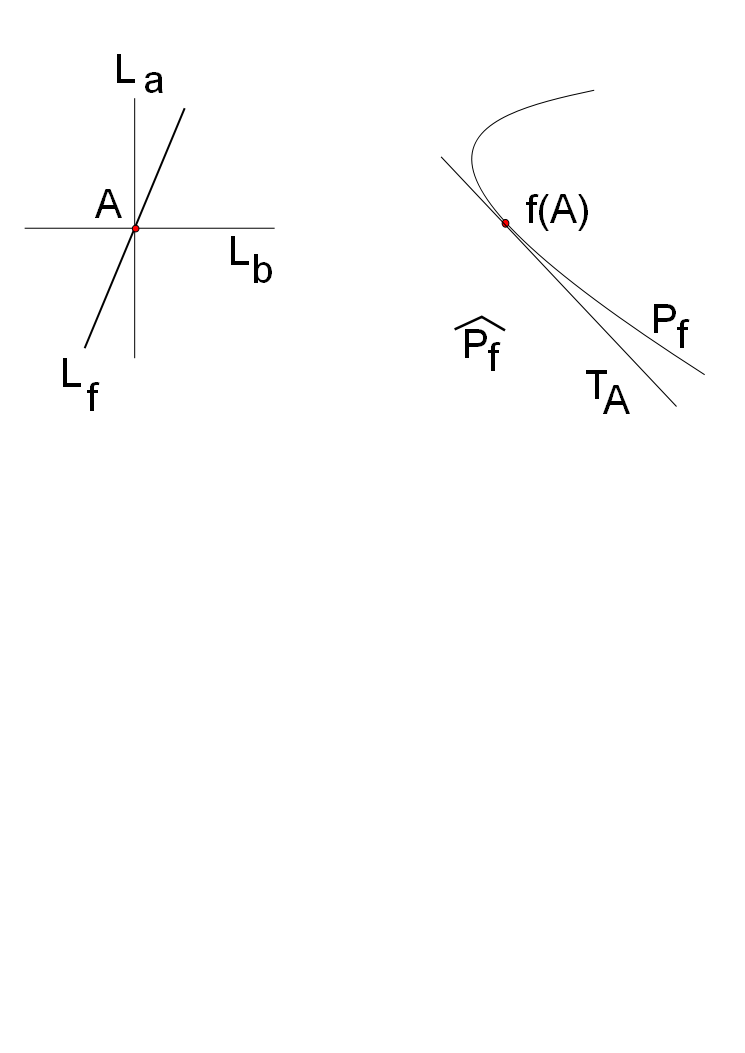}
\vskip -6cm
\caption{Folding parabola $P_f$, parabolic region $\widehat{P_f}$,  tangent line $T_A$.}
\label{3-1}
\end{center}
\end{figure}

\begin{figure}[h!]  
\vskip -2mm
\begin{center}
\includegraphics[width=4in, keepaspectratio] {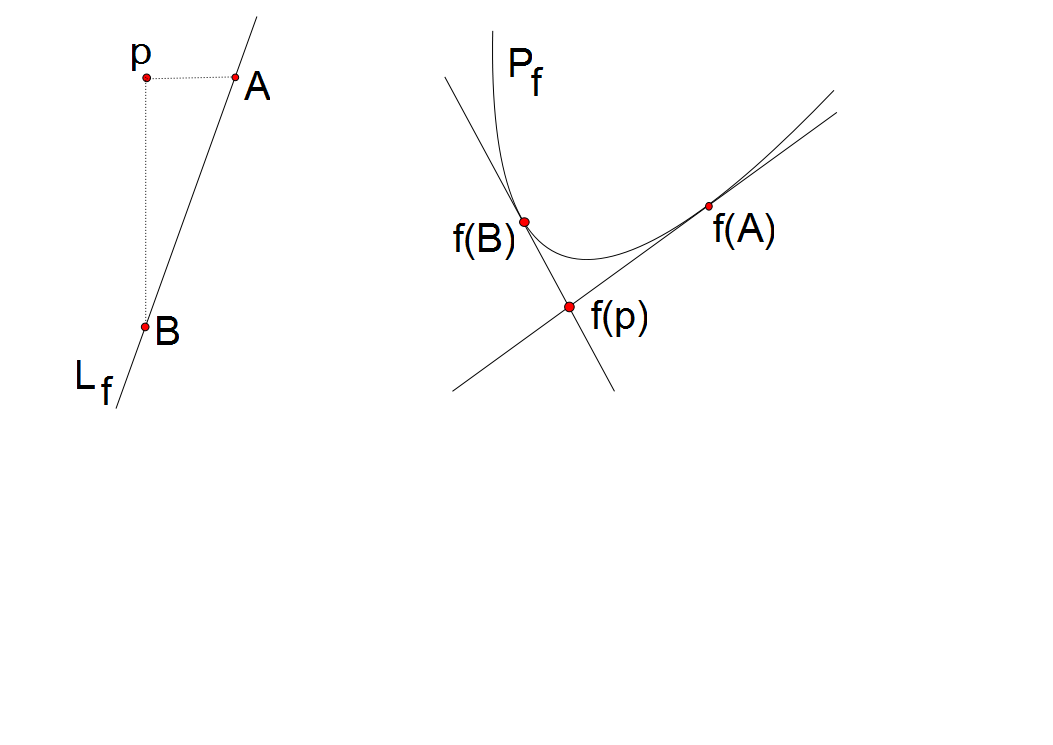}
\vskip -3cm
\caption{The image of a point $\bb p$ under a bi-affine function.  The point $f(\bb p)$ has  ``coordinates"  $(A,B)$.}
\label{4-1}
\end{center}
\end{figure}

\noindent This somewhat complicated formula is merely an analytic expression of the simple geometry shown in Figure~\ref{2-1}.  Statement 3 in Theorem~\ref{image} below makes precise
the 2-to-1 nature of a bi-affine function given in statement 4 of Theorem~\ref{properties}.  
 Statement 1 in Theorem~\ref{image} is illustrated in Figure~\ref{3-1}.  Statement 2 gives a geometric construction
of the image of a given point under a bi-affine function and is illustrated in Figure~\ref{4-1}.  

\begin{theorem} \label{image} Assume that $f$ is a non-degenerate bi-affine map.  We use the notation
$L_a$ for the line $x = a$ and $L_b$ for the line $y=b$.
\begin{enumerate}
\item If $A = (a,b)$ is any point on the folding line $L_f$ and $T_{A}$ is the tangent line
to $P_f$ at $f(A)$, then  $T_A = f(L_{a}) = f(L_{b})$.
\item If $\bb p = (a,b) \in \R^2$, then  $f({\bb p}) =  T_{A} \cap T_{B}$,  where 
$A = L_f \cap L_{a}$ and $B = L_f \cap L_{b}$. 
\item  $f(L_{+}) = f(L_{-}) = \widehat {P_f}$ and, in particular, $f(\bb p) = f({\bb p} ^*)$ for all $\bb p = (x,y) \in \R^2$.
\end{enumerate}
\end{theorem}

\proof
Concerning statement 1, any horizontal or vertical line $L$ intersects $L_f$ in a single point. Therefore $f(L)$ intersects $f(L_f) = P_f$ in a single point, which implies that $f(L)$ is tangent to $P_f$.  Since the intersections $P_f \cap f(L_{a})$ and $P_f \cap f(L_{b}$) both consist of the same single point, $f(L_{a})$ and $f(L_{b})$  both equal the tangent line $T_A$ to $P_f$ at $f(A)$.   

By statement 1, the point $f(\bb p)$ lies both on the tangent to $P_f$ at $f(A)$ and on the tangent to $P_f$ at $f(B)$.  This proves statement 2.  

Consider $\bb p$ and $\bb p^*$ as in Figure~\ref{2-1}.  Statement 2 implies that $f(\bb p) =  T_{A} \cap T_{B} = f({\bb p} ^*)$. In particular $f(L_{+}) = f(L_{-})$.  Since the union of all tangents to
$P_f$ is the parabolic region $\widehat {P_f}$ of $f$, we have $f(L_{+}) = f(L_{-}) = f(\R^2) = \widehat {P_f}$. 
\qed \m

\begin{figure}[!h]  
\begin{center}
\includegraphics[width=4.5in, keepaspectratio] {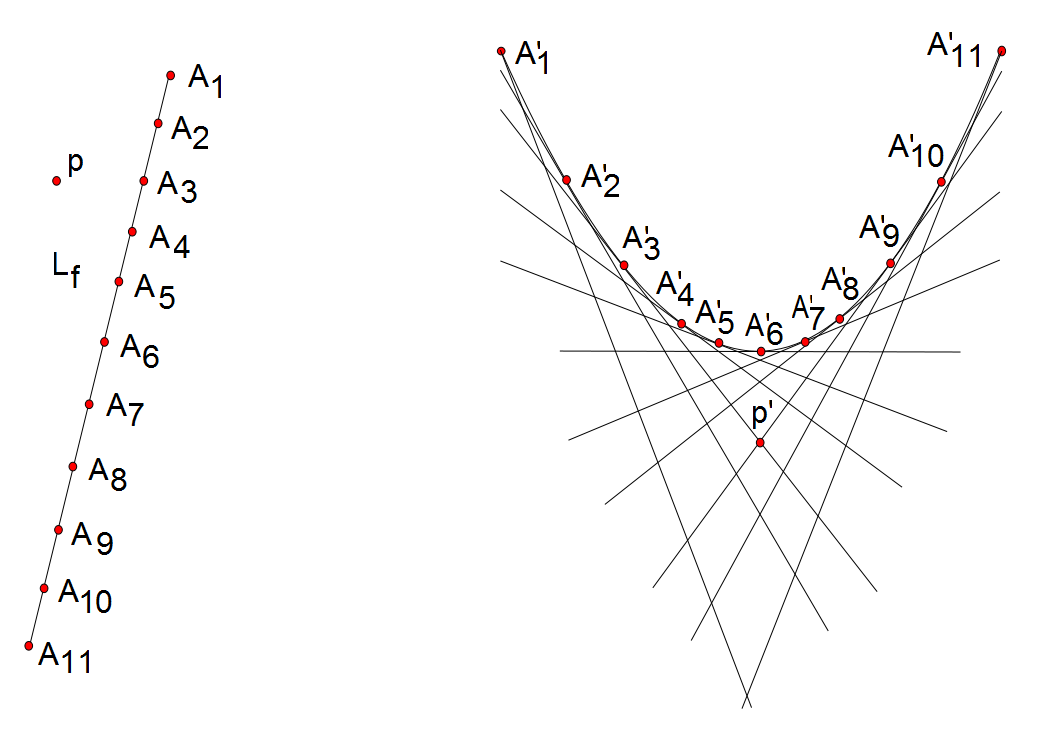}
\caption{Coordinatization of $\widehat{P_f}$.  Note that $A'_i = f(A_i), \, 1\leq i \leq 11,$ and ${\bb p}' = f({\bb p})$.
  The point $\bb p'$ has coordinates $\{A3,A9\}$. }
\label{5-1}
\end{center}
\end{figure}

It is a consequence of Theorem~\ref{image} that the parabolic region can be coordinatized as follows. Each point ${\bb p} \in \widehat{P_f}$ has a unique set $\{A,B\}$ of (unordered) coordinates where $A$ and $B$ are points on the folding line. Specifically $\bb p = T_A \cap T_B$.  This is illustrated in Figures~\ref{4-1} and \ref{5-1}. 

\begin{figure}[htb]  
\begin{center}
\includegraphics[width=4in, keepaspectratio] {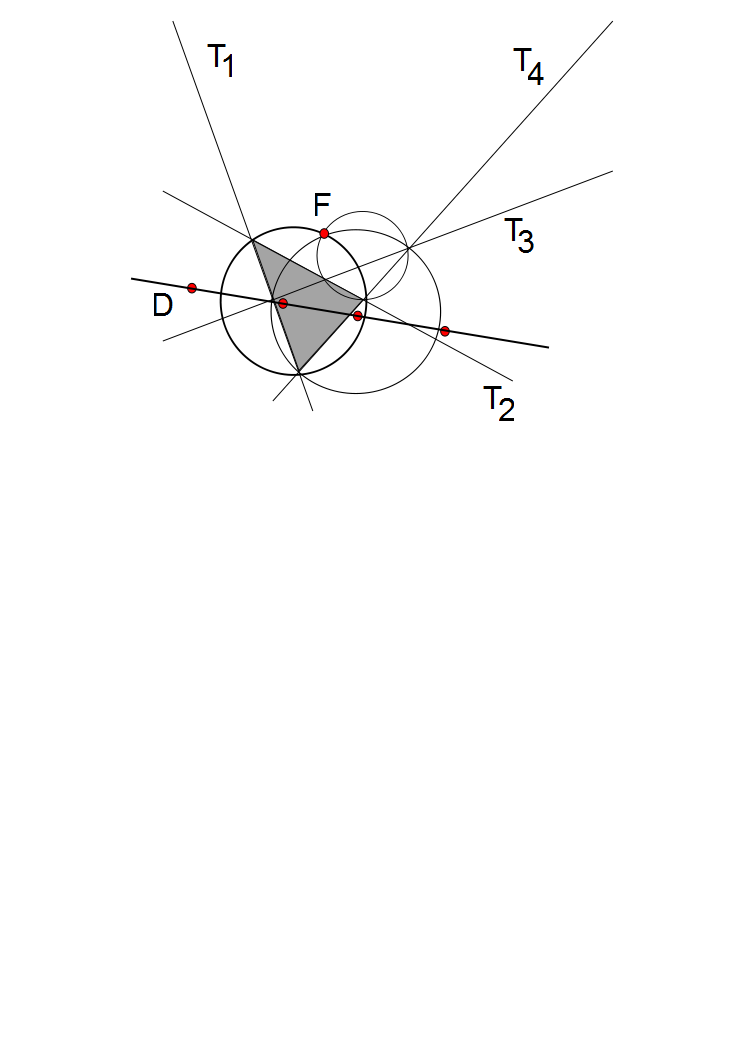}
\vskip -8cm
\caption{Construction of the folding parabola: focus $F$ and directrix $D$.}
\label{6-1}
\end{center}
\end{figure}

The folding parabola $P_f$ can be constructed geometrically as follows (See Figure~\ref{6-1}). Choose four pair wise distinct points $A_1, A_2, A_3, A_4$ on the folding line and construct the four tangent lines $T_{A_1}, T_{A_2}, T_{A_3}, T_{A_4}$.  It is a direct consequence of the classic Brianchon Theorem (actually the converse) that there is a unique parabola with these lines as tangents.  According to what has been shown above, this must be the folding parabola.  The parabola can be explicitly constructed using a theorem of Lambert.  According to Lambert, the circumcircle of a tangent triangle of a parabola (see Figure 6) goes through the focus of the parabola. So the focus $F$ is determined as the intersection of three such circumcircles.  Then reflect $F$ about two of the tangents to get two points on the directrix. 

\section{Iterated Function Systems} \label{IFS}

This section reviews the standard notation and definitions related to iterated function systems (IFS).  These concepts
are then applied to the bi-affine case.

Let $\mathbb{X}$ be a complete metric space. If $f_{m}:\mathbb{X}
\rightarrow\mathbb{X}$, $m=1,2,\dots,M,$ are continuous mappings, then
$\mathcal{F}=\left(  \mathbb{X};f_{1},f_{2},...,f_{M}\right)$ is called an
{\it iterated function system} (IFS). An iterated function system that consists of
bi-affine functions will be called a {\it bi-affine IFS}.  To define the attractor of an IFS, 
first define
\[
\mathcal{F}(B)=\bigcup_{f\in\mathcal{F}}f(B)
\]
for any $B\subset\mathbb{X}$. By slight abuse of terminology we use the same
symbol $\mathcal{F}$ for the IFS, the set of functions in the IFS, and for the
above mapping. For $B\subset\mathbb{X}$, let $\mathcal{F}^{k}(B)$ denote the
$k$-fold composition of $\mathcal{F}$, the union of $f_{i_{1}}\circ f_{i_{2}%
}\circ\cdots\circ f_{i_{k}}(B)$ over all finite words $i_{1}i_{2}\cdots i_{k}$
of length $k.$ Define $\mathcal{F}^{0}(B)=B.$ A nonempty compact set 
$A\subset\mathbb{X}$ is said to be an {\it attractor} of the IFS $\mathcal{F}$ if
\begin{enumerate}
\item $\mathcal{F}(A)=A$ and
\item $\lim_{k\rightarrow\infty}\mathcal{F}^{k}(B)=A,$ for all compact sets
$B\subset \mathbb X$, where the limit is with respect to the Hausdorff metric.
\end{enumerate}
\m

\begin{figure}[htb]
\begin{center}
\includegraphics[width=4cm, keepaspectratio]{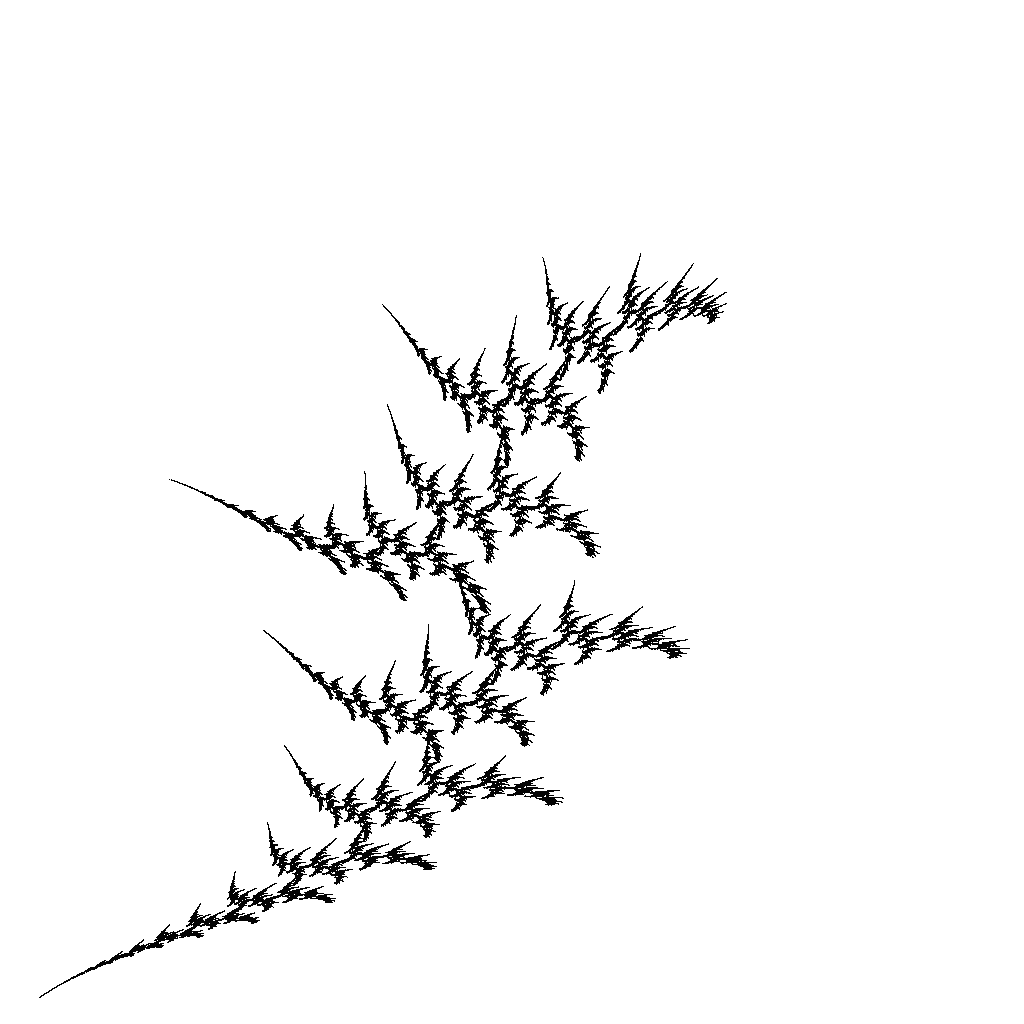} \hskip 8mm \includegraphics[width=4cm, keepaspectratio]{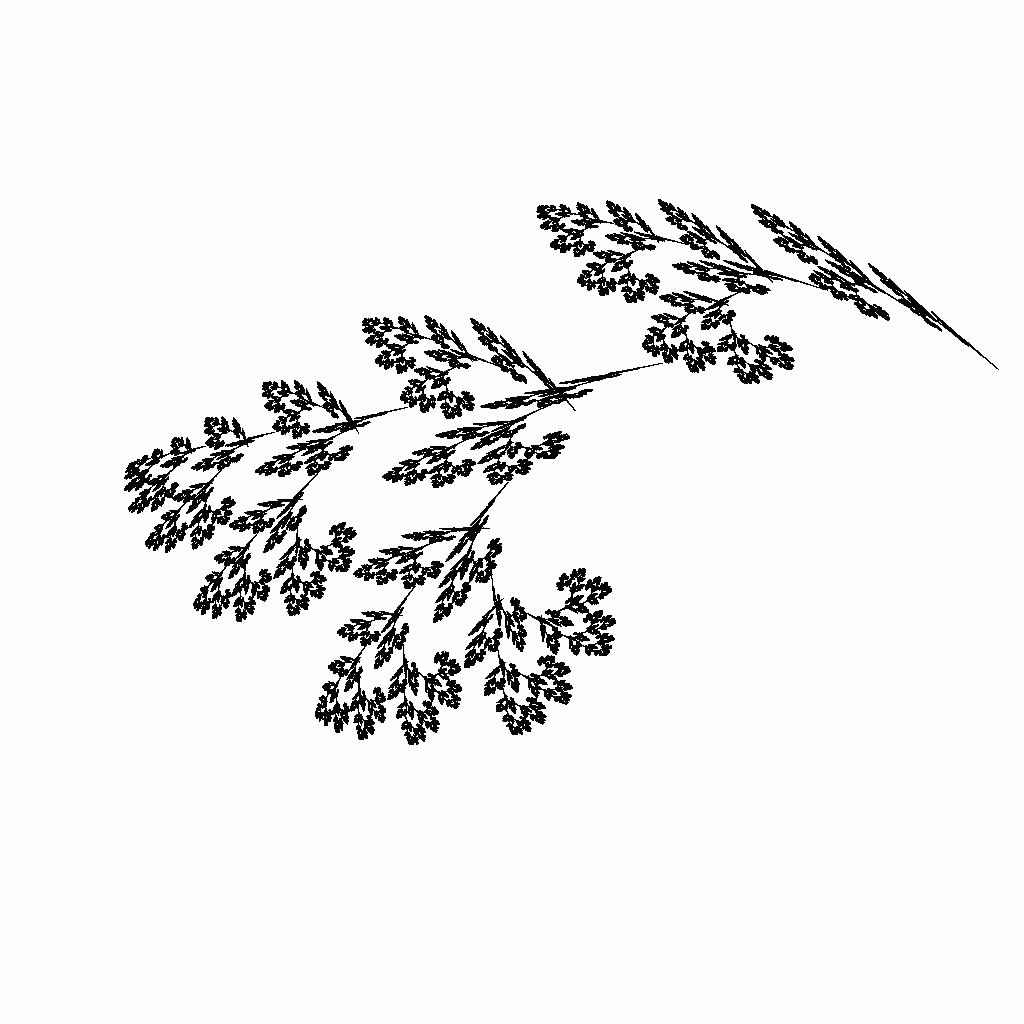}  \hskip 8mm
\includegraphics[width=4cm, keepaspectratio]{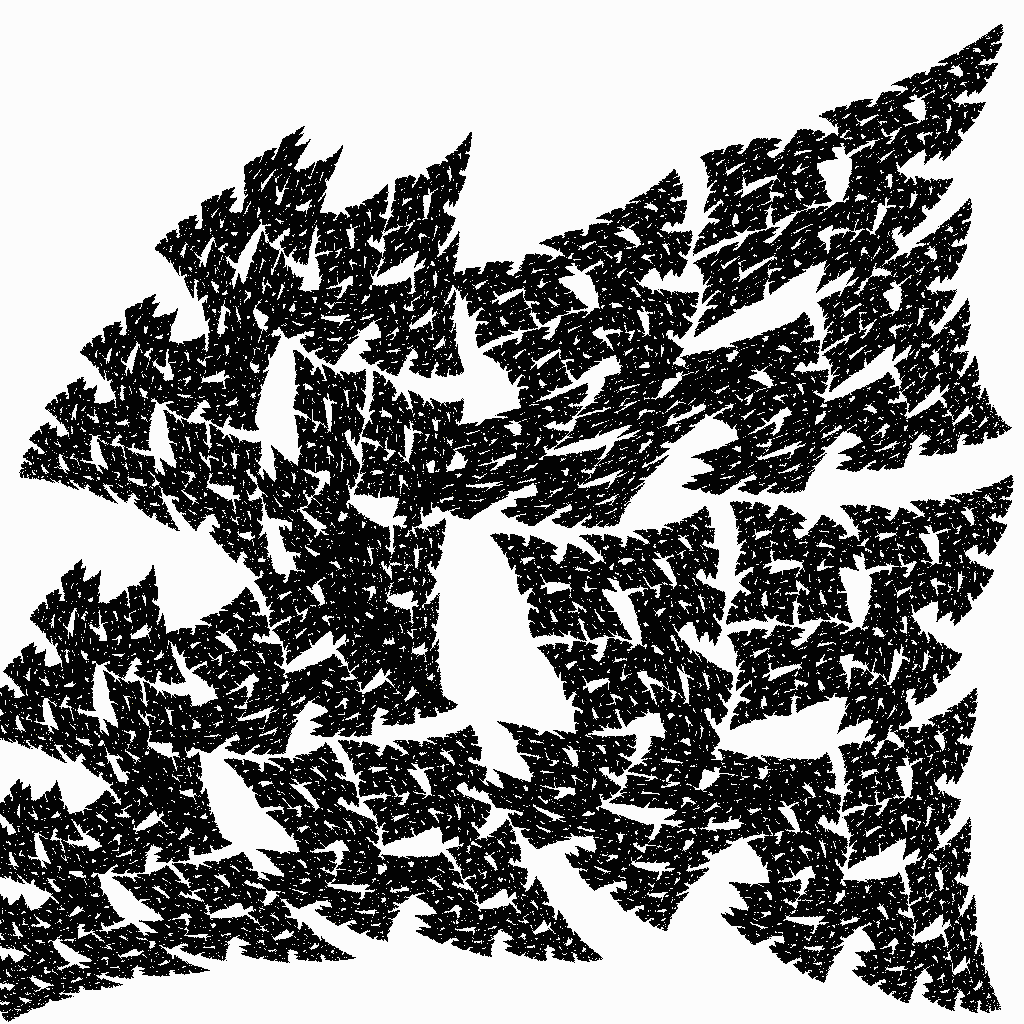} 
\caption{Attractors of  IFSs consisting of two, three, and four bi-affine functions, respectively.}
\label{1-2}
\end{center}
\end{figure}

Attractors for bi-affine IFSs consisting, respectively of $2,3$ and $4$, functions are shown in Figure~\ref{1-2}.

A function $f:\mathbb{X}\rightarrow\mathbb{X}$ is called a
{\it contraction} with respect to a metric $d$ if there is an $0\leq s <1$
such that $d(f(x),f(y))\leq s \,d(x,y)$ for all $x,y\in{\mathbb{R}}^{n}$.  An IFS with the property that each function is a contraction will be called a {\it contractive} IFS.  In his seminal paper  Hutchinson \cite{H} proved  that 
a contractive IFS on a complete metric space has a unique attractor.  Theorem~\ref{contraction} below gives fairly general
conditions under which a bi-affine IFS is contractive.  

\begin{figure}[!h] 
\begin{center}
\vskip -1cm
\includegraphics[width=3in, keepaspectratio] {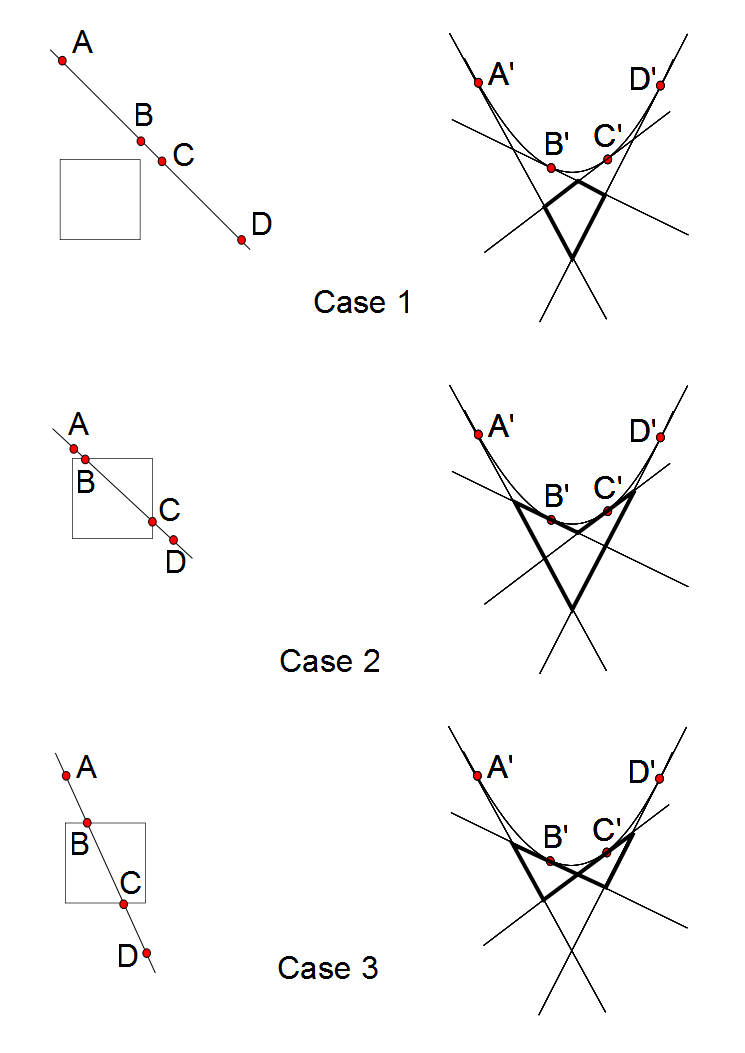}
\caption{The image of the sides of the unit square (on the left) is shown by thick lines (on the right):  $A' = f(A), \,
B' = f(B), \,C' = f(C), \,D' = f(D)$. }
\label{1-3}
\end{center}
\end{figure}

Let $\square$ denote the unit square with vertices $(0,0), (1,0), (1,1), (0,1)$.  The shape of
the image of $\square$ under a non-degenerate bi-affine map $f$ depends on the location of the
folding line $L_f$ relative to $\square$.  It follows from Theorem~\ref{image} that there are
three possible cases as shown in Figure~\ref{1-3}.  It is only in Case 1 ($L_f$ disjoint from the
interior of $\square$) that the image of the sides of $\square$ form a convex quadrilateral.  If this is the case, call
$f$ {\it proper}. \m

Note that the unique bi-affine function taking  $(0,0), (1,0), (1,1), (0,1)$ to the points  $\bb p_0, \bb p_1,\bb p_2,\bb p_3$, respectively, is  
 \begin{equation} \label{quad} f(x,y) = \bb p_0 + (\bb p_1-\bb p_0) x + (\bb p_3-\bb p_0) y 
+ (\bb p_2+ \bb p_0 - \bb p_1 - \bb p_3)xy. \end{equation}
  
\begin{theorem} \label{contraction} Let $f(x,y) = \bb p_0 + (\bb p_1-\bb p_0) x + (\bb p_3-\bb p_0) y + (\bb p_2+ \bb p_0 - \bb p_1 - \bb p_3)xy$
be a proper, non-degenerate bi-affine function. If there is an $s, \, 0 \leq s < 1$, such that 
$$|\bb p_{i+1} - \bb p_i| \leq s, \quad |\bb p_{i+2} - \bb p_{i}| \leq \sqrt{2} s, \quad |\bb p_{i+1} + \bb p_{i-1} - 2 \bb p_i| \leq \sqrt{2} s$$ for $i = 0,1,2,3 \; (mod \; 4)$, then $f$ is a contraction on  $\square$. 
\end{theorem}

In terms of the quadrilateral $\bb p_0 \bb p_1 \bb p_2 \bb p_3$, the first of the three inequalities states that each side has length less than or equal to $s$, the second that each diagonal has length less than or equal to $\sqrt{2} s$, and the third that vector sum of any two incident sides has length less or equal to $\sqrt{2} s$.  For a bi-affine function taking $\square$ into itself, for example, these conditions are not too restrictive.  Two lemmas help in proving Theorem~\ref{contraction}.  

\begin{lemma} \label{lem1} If $f$ is a proper, non-degenerate bi-affine function, then $f$ is injective when restricted to $\square$.  
\end{lemma}

\proof By the comments above, the folding line $L_f$ lies outside the interior of the square $\square$.  The lemma then follows by statement 4 of Theorem~\ref{properties}. 
\qed

\begin{lemma} \label{lem2} Let $f$ be a non-degenerate bi-affine function that is injective on $\square$, and let $W \subseteq \square$ be any rectangle with sides parallel to the $x$ and $y$-axes and with diagonal of length $\rho(W)$.  Further let  $\rho_1(W)$ and $\rho_2(W)$ be the lengths of the two diagonals of $f(W)$ and let 
$$M(\rho_0) = \max_{W : \rho(W) = \rho_0} \max (\rho_1(W),\rho_2(W)).$$
 For any $\rho_0 \leq \sqrt{2}$, if $W_0$ is a rectangle that maximizes $M(\rho_0)$, then $W_0$ and $\square$ have a common vertex.     
\end{lemma}

\proof If one side of $W_0$ lies on the line $x=0$ or $x=1$ and another side lies on $y=0$ or $y=1$, then the proof is complete.  So, without loss of generality, assume that $W_0$ has no side that lies on $x=0$ or $x=1$. Let $V_0$ be the rectangle bounded by the lines $x=0, x=1$ and the lines determined by the upper and lower sides of $W_0$.  Let $A,B,C,D$ be the vertices of $f(V_0)$.  Then by the conditions 1 in the Introduction, there is an $\alpha$ such that the four vertices of $f(W_0)$ are $(1-\alpha)A + \alpha B,  (1-\alpha - \Delta )A + (\alpha + \Delta) B, (1-\alpha- \Delta)C + (\alpha + \Delta) D,(1-\alpha)D + \alpha C$, where $\Delta$ is the horizontal length of $W_0$.  As $\alpha$ varies between $0$ and $\alpha_0$, the rectangle $W_0$ shifts left or right, from the extreme left side of $\square$ to the extreme right side of $\square$.  The lengths of the two diagonals of $f(W_0)$ are $|U_1 + (B+C-A-D)\alpha|$ and $|U_2 + (B+C-A-D)\alpha|$, where vectors $U_1$ and $U_2$ depend on $A,B,C,D$ and $\Delta$. As $\alpha$ varies 
in the range $0 \leq \alpha \leq \alpha_0$, the quantities $U_1 + (B+C-A-D)\alpha$ and $U_2 + (B+C-A-D)\alpha$ describe (parallel) line segments. Hence the maximum of $|U_1 + (B+C-A-D)\alpha|$ and $|U_2 + (B+C-A-D)\alpha|$ occur at an end, i.e. $\alpha = 0$ or $\alpha = \alpha_0$, contradicting the assumption that  $W_0$ has no side that lies on $x=0$ or $x=1$.
\qed \m

\proof (of Theorem~\ref{contraction}) By the continuity of $f$ and the compactness of $\square$, it is sufficient to show that there is an $s, \, 0 \leq s < 1$ such that $|f(x,y) - f(x',y')| \leq s \,|(x,y)-(x',y')|$.  Let $W_0$ be the rectangle whose diagonal is the line segment joining $(x,y)$ and
$(x'y')$. By Lemma~\ref{lem1}, the function $f$ is injective on $\square$, and by Lemma~\ref{lem2}, the 
quantity $$R := \frac{|f(x,y) - f(x',y')|}{|(x,y)-(x',y')|}, \, (x,y) \neq (x',y')$$ is maximized when 
$W_0$ and $\square$ have a common vertex, i.e. when $W_0$ lies on the corner of $\square$.  Without loss of generality, it may be assumed to be the lower left corner. Otherwise, replace $f$ with the composition of $f$ with the rotation that moves the lower left corner to the relevant corner.  Now let $\Delta x = |x'-x|$
and $\Delta y = |y'-y|$.  With $\bb b = \bb p_1 - \bb p_0, \bb c = \bb p_3 - \bb p_0, \bb d = 
\bb p_2 + \bb p_0 - \bb p_1 - \bb p_3$ and setting $r = \frac{\Delta y}{ \Delta x}$, we have two possible
formulas for $R$:  
$$R^2 = \frac{|f(\Delta x, \Delta y) - f(0,0)|^2}{\Delta x, \Delta y)|^2} = 
\frac{1}{1+r^2} | \bb b + r\bb c\ + r \Delta x \bb d |^2 $$ or 
$$R^2 = \frac{|f(\Delta x, 0) - f(0,\Delta y )|^2}{\Delta x, \Delta y)|^2} = \frac{1}{1+r^2} | \bb b - r\bb c |^2.$$  Since $0<\Delta x \leq 1$, the quantity $R^2$ in the first formula is maximized when either
$\Delta x = 0$ or $\Delta x = 1$. So it is now sufficient to prove that 
$$\max_{r>0} \, \big ( \frac{1}{1+r^2} | \bb b + r\bb c|^2, \, \frac{1}{1+r^2} | \bb b + r\bb c + r  \bb d |^2 , \, \frac{1}{1+r^2} | \bb b - r\bb c|^2 \big )  < 1.$$    The inequalities in the statement of the theorem, in terms of $\bb b, \bb c, \bb d$ are that $|\bb b|, |\bb c|, |\bb b+\bb d|, |\bb c+\bb d|$ are all less than or equal to $s$ and
$|\bb b-\bb c|, |\bb b+\bb c+\bb d|, |\bb b+\bb c|,|\bb b-\bb c+\bb d|, |\bb c-\bb b+\bb d|, |\bb b+\bb c+2\bb d|$ are all less than or equal to $\sqrt{2}s$. Maximizing
the three quantities over $r$, in order to verify the above inequality, is a not quite trivial calculus problem whose details are omitted.  
\qed 

\begin{figure}[htp]
\vskip -5mm
\begin{center}
\includegraphics[width=4in, keepaspectratio] {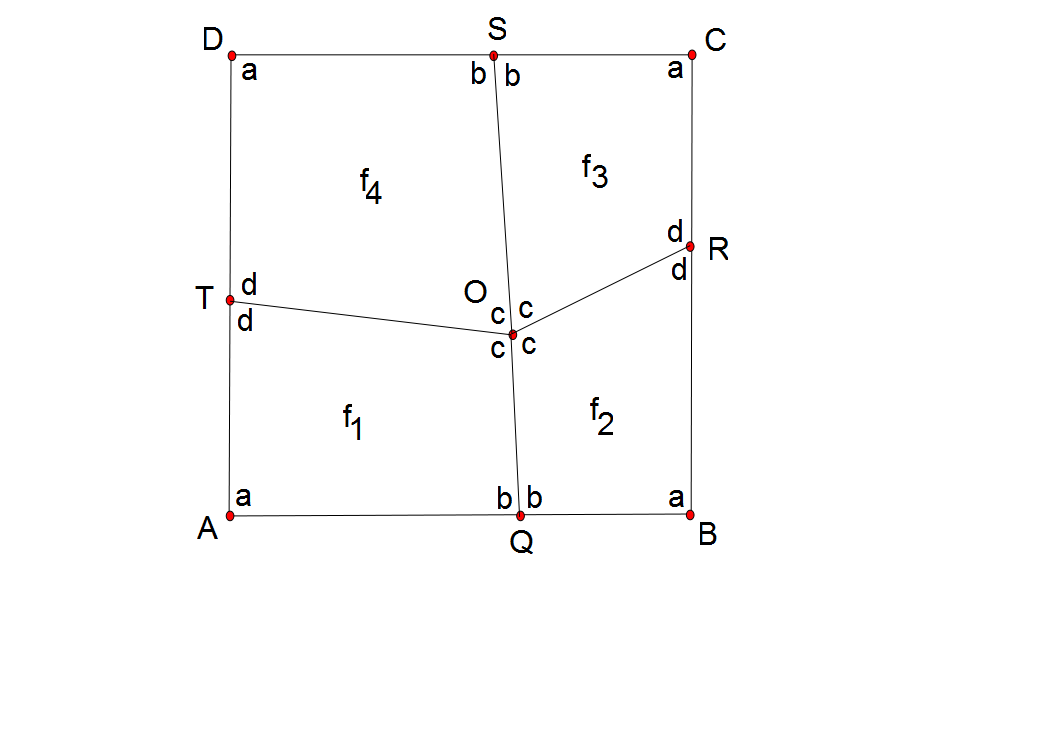}
\vskip -2cm
\caption{An IFS consists of four bi-affine functions. Points labeled with lower case letters are the images of points labeled with upper case letters. The attractor is the unit square. }
\label{2-2}
\end{center}
\end{figure}

\begin{example}  \label{EX1} Consider a bi-affine IFS ${\mathcal F} = \{ \square \, ; \, f_1, f_2, f_3, f_4\}$, where the four functions are determined by the images of the four vertices of $\square$ as shown in Figure~\ref{2-2}.  Each of the four images $f_i(\square),\, i=1,2,3,4$, of the square contains exactly one vertex of $\square$.   The attractor of ${\mathcal F}$ is $\square$ itself.  For ``most" choices of the center and side points, ${\mathcal F}$ satisfies the conditions of Theorem~\ref{contraction} and hence ${\mathcal F}$ is a contractive IFS.   Why this simple example should be of interest is the subject of the next two sections.
\end{example}  

\section{Fractal Homeomorphism} \label{FT}

In the case of a contractive IFS, it is possible to assign to each point of the attractor an ``address".  Given two IFSs ${\mathcal F}$ and ${\mathcal G}$ with respective attractors $A_F$ and $A_G$, a fractal homeomorphism is basically a homeomorphism $h \, : \, A_F \rightarrow A_G$ that sends a point in $A_F$ to the point in $A_G$ with the same address. 
To make this notion  precise, let $\mathcal{F}=\left(  \mathbb{X};f_{1},f_{2},...,f_{N}\right)$  be a contractive IFSs on a complete metric space $\mathbb{X}$ with attractor $A$.  Let 
$\Omega = \{1,2, \dots ,N\}^{\infty}$ denote the set of infinite strings using symbols $1,2, \dots , N$, and for
$\sigma \in \Omega$, let $\sigma|_k$ denote the string consisting of the first $k$ symbols in $\sigma$.  Moreover,
if $\sigma = i_0\, i_1, \, i_2 \cdots$, then we use the notation $$f_{\sigma|_n} := f_{i_0} \circ f_{i_1}\circ \cdots \circ f_{i_n}.$$
Now define a map $\pi: \Omega \rightarrow A$, called the {\it coding map}, by
 $$\pi(\sigma) :=
\lim_{k\rightarrow\infty}f_{\sigma|_k}(x).$$ It is well known \cite{A,H} that the limit above exists and is independent of $x \in \mathbb X$.  Moreover $\pi$  is continuous, onto, and satisfies the following commuting diagram  for each $n = 1,2,\dots, N$.
\begin{equation} \label{codingDiagram}
\begin{array}
[c]{ccc}
\Omega & \overset{s_{n}}{\rightarrow} & \Omega\\
\pi\downarrow\text{\ \ \ \ } &  & \text{ \ \ \ }\downarrow\pi\\
\mathbb X & \underset{f_{n}}{\rightarrow} & \mathbb X
\end{array}
\end{equation}
The symbol $s_{n}:\Omega\rightarrow\Omega$ denotes the inverse shift map defined by $s_{n}(\sigma)=n \sigma.$
 A {\it section} of the coding map $\pi$ is a function $\tau \, : \, \Omega \rightarrow A$ 
such that $\pi \circ \tau$ is the identity. A section selects, for each $x \in A$, a single
{\it address} in $\Omega$ from the ones that come from the coding map.  Call the set $\Omega_{\tau} := \tau(A)$ the {\it address space} of the section $\tau$.   Let $S$ denote the shift operator on $\Omega$, i.e, $S(n \sigma) = \sigma$ for any $n \in \{ 1,2, \dots , N\}$ and any $\sigma \in \Omega$.  A subset $W \subset \Omega$ will be called {\it shift invariant} if $\sigma \in W$ implies that $S(\sigma) \in W$.  If $\Omega_{\tau}$ is shift invariant, then $\tau$ is called a {\it shift invariant section}.  The following example demonstrates the naturalness of shift invariant sections. 

\begin{example} \label{EX0}  Consider the IFS ${\mathcal F} = ( \R \, ; \, f_0,f_1)$ where $f_0(x) = \frac12\, x$ and $f_1(x) = \frac12\, x + \frac12$.  The attractor is the interval $[0,1]$.  An address of a point $x$ is a binary representation of $x$.  In choosing a section $\tau$ one must decide, for example, whether to take $\tau (\frac14) = .01$ or $\tau (\frac14) = .00111\cdots$. If the section $\tau$ is shift invariant, this would imply, for example, that if $\tau(\frac14) = .00111\cdots$, then  $\tau (\frac12) = .0111\cdots$, not  $\tau (\frac12) = .100\cdots$.
\end{example} 

\begin{lemma} \label{shift} With notation as above,  a section $\tau$ of an IFS is shift invariant if and only if, for any $x\in A$, if $\tau(x)|_1 = n$, then $(S\circ \tau)(x) = (\tau\circ f^{-1}_n)(x)$.  
\end{lemma}

\proof Given the right hand statement above,  we will prove that $\tau$ is shift invariant.  Assume that $n\sigma \in \Omega_{\tau}$.  Then there is an $x\in A$ such that $\tau(x) = n\sigma$, and hence $\sigma = (S\circ \tau)(x) =
 (\tau\circ f^{-1}_n)(x) = \tau ( f^{-1}_n)(x))$.  Thus $\sigma \in \Omega_{\tau}$.

Conversely, assume that $\tau$ is shift invariant.  Assume that $x\in A$ and $\tau(x) = n \sigma$ for some $\sigma \in \Omega$.  By shift invariance, there is a $y \in A$ such that $\tau(y) = \sigma$.  Now
$$x = (\pi\circ\tau)(x) = \pi (n \sigma) = f_n (\lim_{k \rightarrow \infty} f_{\sigma|_k}(A)) = f_n(\pi(\tau(y))) = f_n(y).$$
Therefore $y = f_n^{-1}(x)$ and $$(S\circ \tau)(x) = S(n\sigma) = \sigma = \tau(y) = \tau(r_n^{-1}(x)) = (\tau\circ f^{-1}_n)(x).$$  
\qed 

Call an IFS {\it injective} if each function in the IFS is injective.  Theorem~\ref{mask} below states that every shift invariant section of an injective IFS can be obtained from a mask.
For an IFS $\mathcal F$ with attractor $A$, a {\it mask} is a partition $M = \{M_i, 1\leq i \leq N\}$ of $A$ such that $M_i \subseteq f_i(A)$ for all $f_i\in {\mathcal F}$.    Given an injective IFS $\mathcal F$ and a mask $M$, consider the function $T \, : \, A \rightarrow A$ defined by $T(x) := f_i^{-1} (x)$ when $x\in M_i$.  The {\it itinerary} $\tau_M(x)$ of a point $x\in A$ is the string $i_0\, i_1 \, i_2 \cdots \in \Omega$, where $i_k$ is the unique integer $1\leq i_k \leq N$ such that $$T^k(x) \in M_{i_k}.$$  

\begin{theorem} \label{mask} Let $\mathcal F$ be a contractive and injective IFS. 
\begin{enumerate}
\item If $M$ is a mask, then $\tau_M$ is a shift invariant section of $\pi$. 
\item If $\tau$ is a shift invariant section of $\pi$, then $\tau = \tau_M$ for some mask $M$.  
\end{enumerate} 
\end{theorem} 

\proof To show that $\tau_M$ is a section, let $x\in A$.  Then 
$(\pi\circ \tau_M)(x) = \lim_{k\rightarrow \infty} f_{\tau_M(x)|_k}(A)$. It follows immediately from the definition of $\tau_M$  that $x \in f_{\tau_M (x)|_k}(M_{i_{k+1}}) \subseteq  f_{\tau_M (x)|_k}(A) $ for all $k$. Hence $(\pi\circ \tau_M)(x) = x$.
Concerning the shift invariance, it follows from the definition of $T$ that the following diagram commutes.
\begin{equation} 
\begin{array}
[c]{ccc}
A & \overset{T}{\rightarrow} & A \\
\tau_M\downarrow\text{\ \ \ \ } &  & \text{ \ \ \ }\downarrow\tau_M \\
\Omega & \underset{S}{\rightarrow} & \Omega
\end{array}
\end{equation}
If $\sigma \in \Omega_{\tau_M}$, then there is an $x \in A$ such that $\sigma = \tau_M(x)$ and, from the diagram, 
$S(\sigma) = (S\circ \tau_M)(x) = (\tau_M \circ T)(x) =  \tau_M (f_n^{-1} (x)) \in \Omega_{\tau_M}$ for some $n \in \{ 1,2, \dots, N\}$.

Concerning the second statement, define a mask $M = \{M_i, 1\leq i \leq N\}$ as follows:
$$M_i = \{ x \, : \, \tau (x) = i \, \sigma \quad \text{for some} \quad  \sigma \in \Omega \}. $$
  It is sufficient to show that $M_i \subseteq f_i(A)$, and that $\tau_M(x) = \tau(x)$ for all $x\in A$. If $x\in M_i$, then $\tau(x) = i\sigma$ for some $\sigma \in \Omega$ and $x = (\pi \circ \tau )(x) = f_i(\lim_{k\rightarrow \infty} f_{\sigma|_k} (A)) \in f_i(A)$.  

To show that $\tau_M(x) = \tau(x)$, let $\tau(x) = j_0 j_1 j_2 \cdots$  and $\tau_M(x) = k_0 k_1 k_2 \cdots$.   That $j_0  = k_0$  follows form the definitions.  By induction, assume that $j_i = k_i\, i=0,1,\dots , m-1$.   Applying Lemma~\ref{shift} for $m$ times yields $$j_m =( S^m\circ \tau)(x) )|_1 = \tau (f_{j_{m-1}}^{-1} \circ \cdots \circ f_{j_1}^{-1}
\circ f_{j_0}^{-1}(x) |_1,$$  
where the $j_i$'s are determined by the recursive formula  $\tau (f_{j_{r-1}}^{-1} \circ \cdots \circ f_{j_1}^{-1}
\circ f_{j_0}^{-1}(x) |_1 = j_r$. 
By the definition of the mask, $ \tau (f_{j_{m-1}}^{-1} \circ \cdots \circ f_{j_1}^{-1}
\circ f_{k_0}^{-1}(x) ) \in M_{j_m}$.  But by the definition of the itinerary, the $k_i$'s are determined by the recursive formula 
$$ \tau (f_{k_{r-1}}^{-1} \circ \cdots \circ f_{k_1}^{-1}\circ f_{k_0}^{-1}(x))  \in M_{k_{r}}$$
for $r = 0,1,2, \dots , m$.  But, since  $k_i = j_i\, i=0,1,\dots , m-1$, we have $ \tau (f_{j_{m-1}}^{-1} \circ \cdots \circ f_{j_1}^{-1}\circ f_{j_0}^{-1}(x) ) \in M_{k_m}.$  Therefore $k_m = j_m$.
\qed \m 

To define fractal homeomorphism, consider two contractive IFSs $\mathcal{F}=\left(  \mathbb{X};f_{1},f_{2},...,f_{N}\right)$ and $\mathcal{G}=\left(  \mathbb{X};g_{1},g_{2},...,g_{N}\right)$  with the same number $N$ of functions on a complete metric space
$\mathbb{X}$.  Let $A_F$ and $A_G$ be the  attractors and $\pi_F$ and $\pi_G$ the coding maps of $\mathcal F$ and $\mathcal G$, respectively.  A homeomorphism $h \, : \, A_F \rightarrow A_G$ is called a {\it fractal homeomorphism } if
there exist shift invariant sections $\tau_F$ and $\tau_G$ such that the following diagram commutes: 
\begin{equation} \label{homeo}
\begin{array}
[c]{ccc}
A_F & \underset{h}{\rightarrow} & A_G \\ \text{$\tau_F$} \searrow & & \swarrow \text{$\tau_G$} \\
 & \Omega & 
\end{array}
\end{equation}
i.e., the homeomorphism $h$ takes each point $x\in A_F$ with address $\sigma = \tau_F(x)$ to the point $y \in A_G$ with the same address $\sigma =\tau_G(y) $.  Theorem~\ref{pi-tau} below states that the fractal homeomorphisms between
attractors $A_F$ and $A_G$ are exactly mappings of the form $\pi_G\circ \tau_F$ or $\pi_F\circ\tau_G$ for some shift invariant sections $\tau_F, \, \tau_G$.   

\begin{theorem}  \label{pi-tau} Let $\mathcal F$ and $\mathcal G$ be contractive IFSs.  With notation as above:
\begin{enumerate} 
\item If $h \, : \, A_F \rightarrow A_G$ is a fractal homeomorphism with corresponding sections $\tau_F$ and $\tau_G$, then $\Omega_{\tau_F} = \Omega_{\tau_G}$.  Moreover $h = \pi_G \circ \tau_F$ and $h^{-1} = \pi_F \circ \tau_G$ .
\item If $\tau_F$ is a shift invariant section for $\mathcal F$ and $h := \pi_G \circ \tau_F$ is a homeomorphism, then $h$ is a fractal homeomorphism.
\end{enumerate}
\end{theorem}

\proof Concerning statement 1, since $h$ is a bijection, the commuting diagram~\ref{homeo} implies that the images of $\tau_F$ and $\tau_G$ are equal, i.e.,  $\Omega_{\tau_F} = \Omega_{\tau_G}$. Now $\tau_F = \tau_G\circ h$ from the diagram implies $\pi_G \circ \tau_F = \pi_G \circ \tau_F =
(\pi_G \circ \tau_G) \circ h = h$.  The formula involving $h^{-1}$ is likewise proved.

Concerning statement 2,  the section $\tau_F$ is a bijection from $A_F$ onto $\Omega_{\tau_F}$.  Since $h$ is also a bijection, the equality  $h = \pi_G \circ \tau_F$ implies that $\pi_G|_{ \Omega_F}$, the restriction of $\pi_G$ to $\Omega_F$, is a bijection onto $A_G$.  If $\tau_G$ is the inverse of $\pi_G|_{ \Omega_F}$, then $\tau_F$ and $\tau_G$ satisfy the commuting diagram~\ref{homeo}.  That $\tau_F$ is shift invariant means that $\Omega_{\tau_G} = \Omega_{\tau_F}$ , i.e., $\tau_G$ is shift invariant.  
\qed \m

\section{Image from a Fractal Homeomorphism} \label{frSection}

 This section concerns images on the unit square $\square$.  Define an {\it image} as
 a function $c \, : \, \square \rightarrow {\cal C}$, where $\cal C$ denotes the color palate, for example
 ${\cal C} = \{ 0,1,2,\dots, 255\}^3$.  If $h$ is any homeomorphism from $\square$ onto $\square$, define the {\it transformed image} $h(c) \, : \, \square \rightarrow {\cal C}$ by
 $$h(c) := c \circ h.$$ 
We are interested in the case where $h$ is a fractal homeomorphism.  The remainder of this section concerns
fractal homeomorphism based on bi-affine IFSs with four functions as described in Example~\ref{EX1}. \m

 Consider  Example~\ref{EX1} depicted in Figure~\ref{2-2}.  For the bi-affine IFS \linebreak
 ${\mathcal F} = \{ \square \, ; \, f_1, f_2, f_3, f_4\}$,  we will construct a section $\tau_F$ that is referred to in \cite{B1} as the {\it top section}.  Consider the mask $M_F = \{M_1, M_2, M_3, M_4\}$ 
defined recursively by $$M_i = f_i(\square)\setminus \bigcup_{j=1}^{i-1} f_j(\square)$$ for $i = 1,2,3,4$. 
Explicitly, $M_1$ is the closed quadrilateral $ATOQ$, $M_2$ is the open quadrilateral $OQBR$ together with the segments $(Q,B],[B,R], [R,O)$,  $M_3$ is the open quadrilateral $ORCS$ together with the segments $(R,C], [C,S],[S,O)$, and
$M_4$ is the open quadrilateral $OSDT$ together with the segments $(S,D],[D,T)$.  The section $\tau_F$ corresponding to the mask $M_F$ is given by $\tau(x) = \max \pi^{-1}(x)$, where the maximum is with respect to the lexicographic order on $\Omega$. 

 Now consider a second bi-affine IFS
${\mathcal G} = \{ \square \, ; \, g_1, g_2, g_3, g_4\}$ of the same type with points $O',Q',R',S',T'$ replacing 
$O,Q,R,S,T$, and with mask $M_G$ defined exactly as it was for $M_F$.  The masks $M_F$ and $M_G$ induce
shift invariant sections $\tau_F$ and $\tau_G$, respectively, as verified by Theorem~\ref{mask}.  Theorem~\ref{picture} below  states  that $\pi_G \circ \tau_F$ and $\pi_F\circ \tau_G$ are  continuous and hence, by Theorem~\ref{pi-tau},  fractal homeomorphisms. 

To prove Theorem~\ref{picture}, the following lemma will be used. The proof is routine and will be omitted.  All  partitions $\mathcal P$ will be of  the unit square $\square$ into regions whose closures are topological polygons.  The {\it dual graph} of such a partition is the graph $\Gamma_{\mathcal P}$ whose points are the regions and where two vertices are joined if and only if the corresponding regions share a side.  A partition $\mathcal P$ is {\it nested} in partition $\mathcal Q$  if each region in $\mathcal P$ is contained in some region of $\mathcal Q$.  Assume that partition $\mathcal {P}_1$ is nested in partition $\mathcal{P}_2$ and $\mathcal{Q}_1$ is nested in $\mathcal{Q}_2$, and that there are graph isomorphisms $\Phi_1 \, :\, \Gamma_{\mathcal{P}_1} \rightarrow \Gamma_{\mathcal{Q}_1}$ and  $\Phi_2 \, :\, \Gamma_{\mathcal{P}_2} \rightarrow \Gamma_{\mathcal{Q}_2}$.  Call $\Phi_1$ and $\Phi_2$ {\it compatible} if whenever  $P_1 \in {\mathcal P}_1$ and 
$P_2 \in {\mathcal P}_2$ with $P_1 \subseteq P_2$ we have $\Phi_1(P_1) \subset \Phi_2(P_2)$.
The {\it mesh}  $|\mathcal P|$ of a partition $\mathcal P$ is the maximum diameter of the regions.  If  
$\lim_{n \rightarrow \infty} |\mathcal P_n| = 0$, then, for any $x \in \square$, there is a unique nested sequence $\{P_n\}$ of regions $P_n \in \mathcal P_n$ such that  $x = \bigcap_{n \in \mathbb N}  P_n$.  

\begin{lemma} \label{lem}  Let $\mathcal P_n$ and  $\mathcal Q_n, \, n = 0,1,2, \dots$,  be two nested sequences of partitions of the unit square $\square$ with $\lim_{n \rightarrow \infty} |\mathcal P_n| = \lim_{n \rightarrow \infty} |\mathcal Q_n| = 0$.   Assume that there are corresponding sequences of compatible graph isomorphisms $\Phi_n \, : \, G_{ \mathcal P_ n} \rightarrow G_{ \mathcal Q_ n}$.   The map $h \, : \, \square \rightarrow \square$ defined as follows is a homeomorphism.   For $x \in \square$, let $x = \bigcap_{n \in \mathbb N}  P_n \in \square$ with $P_n \in \mathcal P_n$, and define $h(x) := \bigcap_{n \in \mathbb N}  \Phi_n(P_n)$.
\end{lemma}

\begin{figure}[htb]
\vskip -.5cm
\begin{center}
\includegraphics[width=2.5in, keepaspectratio] {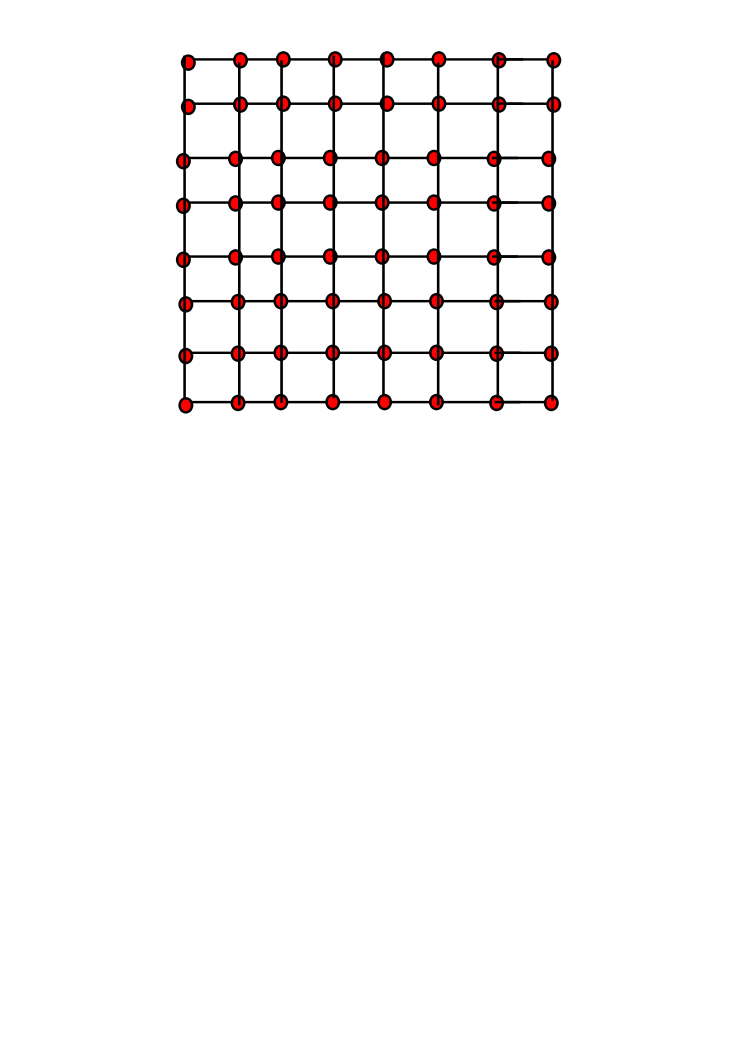}
\vskip -5cm
\caption{The dual graph of the partition $\mathcal P_2$. }
\label{grid}
\end{center}
\end{figure}

\begin{theorem}  \label{picture}
For the two  bi-affine IFS's ${\mathcal F} = \{ \square \, ; \, f_1, f_2, f_3, f_4\}$ and \linebreak
${\mathcal G} = \{ \square \, ; \, g_1, g_2, g_3, g_4\}$ defined above, the map $h = \pi_G \circ \tau_F$ is a homeomorphism.  
\end{theorem}
 
\proof    For each $n\geq 0$, let $\Omega_n$ denote the set of strings of length $n$ using symbols $\{1,2,3,4\}$ For the IFS $\mathcal F$, define a partition ${\mathcal P}_F^n = \{ P_{\sigma} \, : \, \sigma \in \Omega_n\}$ of $\square$ recursively by taking ${\mathcal P}_F^0 = {\mathcal P}_F$ and $${\mathcal P}_F^{n+1} =\{ P_{\sigma \, j} = P_{\sigma} \cap f_{\sigma}(P_j) \, : \, \sigma \in \Omega_n, 1\leq j \leq 4 \}.$$ A straightforward induction shows that $\{ M_F^n\}$ is a nested sequence of partitions of $\square$.  The dual graph $\Gamma_F^n$ of ${\mathcal P}_F^n$ is the grid graph shown in Figure~\ref{grid} for $n = 2$.   This construction of a nested partition can be repeated for the IFS $\mathcal G$.  Since the obvious graph isomorphisms between $\Gamma_F^n$ and $\Gamma_G^n$ are compatible with the nested partitions, Lemma~\ref{lem} implies that $h$ is a homeomorphism.  
\qed

\end{document}